\magnification\magstep 1
\font\twelvebf=cmbx12
\def\Cdb{ \; {}^{ {}_\vert }\!\!\!{\rm C}   }

\def\Qdb{ {\rm Q}\kern-.65em {}^{ {}_/ }}

${}$ \vskip 1cm
\centerline{\twelvebf ON FUNCTION AND OPERATOR MODULES}
\vskip 1cm
\centerline{\bf David BLECHER and Christian LE MERDY}
\vskip .5cm
\centerline{5/18/99}
\vskip 1.5cm
\noindent
{\bf Abstract.} Let $A$ be a unital Banach algebra. We give a characterization
of the left Banach $A$-modules $X$ for which there exists a commutative unital
$C^*$-algebra $C(K)$, a linear isometry $i\colon X\to C(K)$, and a
contractive unital homomorphism $\theta\colon A\to C(K)$ such that
$i(a\cdotp x) =\theta(a)i(x)$ for any $a\in A,\ x\in X$. We then deduce a
``commutative" version of the Christensen-Effros-Sinclair characterization
of operator bimodules. In the last section of the paper, we prove a
$w^*$-version of the latter characterization, which generalizes some
previous work of Effros and Ruan.
\bigskip\medskip\noindent
{\bf 1 - \underbar{Introduction.} }
\medskip\noindent
Let $H$ be a Hilbert space and let $B(H)$ be the $C^*$-algebra of all bounded operators
on $H$. Let $A\subset B(H)$ and $B\subset B(H)$ be two unital closed subalgebras and let 
$X\subset B(H)$ be a closed subspace. If $axb$ belongs to $X$ whenever $a\in A$,
$x\in X$, $b\in B$, then $X$ is called a (concrete) operator $A$-$B$-bimodule. 
The starting point of this paper is the abstract characterization of 
these bimodules due to Christensen, Effros, and Sinclair.
Namely, let us consider two unital operator
algebras $A$ and $B$ and let $X$ be an arbitrary operator space. Assume that
$X$ is an $A$-$B$-bimodule. Then for any integer $n\geq 1$, the Banach space 
$M_n(X)$ of $n\times n$ matrices with entries in $X$ can be naturally
regarded as an $M_n(A)$-$M_n(B)$-bimodule, by letting
$$
[a_{ik}]\cdotp [x_{kl}]\cdotp [b_{lj}] =\Bigl[ \,\sum_{1\leq k,l\leq n} a_{ik}\cdotp 
x_{kl}\cdotp b_{lj}\,\Bigr]
$$
for any $[a_{ik}]\in M_n(A),\ [x_{kl}]\in M_n(X),\ [b_{lj}]\in M_n(B)$. It is shown in [CES]
that if $\Vert a\cdotp x\cdotp b\Vert\leq \Vert a\Vert\Vert x\Vert\Vert b\Vert$ for any
$n\geq 1$ and any $a\in M_n(A),\ x\in M_n(X),\,$ and $b\in M_n(B)$,
then there exist a Hilbert space $H$, and three complete isometries 
$$
J\colon X\to B(H),\ \pi_1\colon A\to B(H),\ \pi_2\colon B\to B(H),
\leqno(1.1)
$$ 
such that $\pi_1, \pi_2$ are homomorphisms and $J(a\cdotp x\cdotp b)
=\pi_1(a)J(x)\pi_2(b)$ for any $a\in A$, $x\in X$, $b\in B$. In that case, 
$X$ is called an (abstract) operator $A$-$B$-bimodule.
Equivalently (in the Operator Space language), 
operator $A$-$B$-bimodules $X$ are characterized by the property that
the bimodule action $\, A\times X\times B\longrightarrow X\,$
extends to a completely contractive map from the Haagerup tensor product
$A\otimes_h X\otimes_h B$ into $X$.
\smallskip\noindent
Let us focus on minimal operator spaces, that is operator spaces
which can be embedded completely isometrically into a commutative $C^*$-algebra.
First note that if $A$ is a unital
operator algebra whose operator space structure is minimal, then it can be
represented completely isometrically and algebraically as a unital subalgebra of
a commutative unital $C^*$-algebra ([Bl1, Theorem 5]). In that case, $A$ will be called a 
{\it function algebra}. Now if $A, B$ are two function algebras, if $X$ is a minimal
operator space, and if $X$ is an operator $A$-$B$-bimodule, it is natural to wonder whether
the corresponding bimodule action can be represented within a commutative $C^*$-algebra.
We will show in Corollary 2.10 below that this is indeed the case. In other words, the 
three complete isometries in (1.1) can be chosen to take their values
in a common unital commutative $C^*$-subalgebra of $B(H)$.
\bigskip\noindent
In Section 2 below, in which we prove the latter result, 
we shall essentially work in a Banach Space framework. 
All Banach spaces considered are complex and are assumed to be
different from $(0)$.
Given a unital Banach algebra $A$ (with a normalized unit denoted by $1$)
and a Banach space $X$, we say
that $X$ is a left Banach $A$-module if $X$ is equipped with a left
module action $A\otimes X\to X,\ (a,x)\mapsto a\cdotp x$, such that $1\cdot x=x$
and $\Vert a\cdotp x\Vert\leq\Vert a\Vert\Vert x\Vert$ for any $x\in X,\, a\in A$.
Right Banach $B$-modules and Banach $A$-$B$ bimodules over a pair of unital Banach
algebras $(A,B)$ are defined similarly. As noticed in [ER1, Section 1], no representation
theorem exists for those Banach modules. In Theorem 2.2, we 
establish a general representation result which partly fills this gap. 
Namely, let $X$ be a left Banach $A$-module and let $m\colon
A\otimes X\to X$ be the corresponding module action. We will show that the contractivity
of $m$ with respect to various tensor norms on $A\otimes X$ is equivalent to the existence
of an isometric embedding  $i\colon X\to C(K)$ of $X$ into a commutative unital
$C^*$-algebra $C(K)$ and of a contractive unital homomorphism $\theta\colon A\to C(K)$
such that $i(a\cdotp x) =\theta(a)i(x)$ for any $a\in A,\ x\in X$.
This result will emphasize the role of the so-called multiplier algebra 
${\rm Mult}(X)$ of a Banach space $X$,
as defined e.g. in [B1, p. 54-55]. We will derive some simple consequences of known
results concerning this algebra.
\smallskip\noindent
In Section 3, we continue our investigation of operator bimodules by showing that
the Effros-Ruan representation theorem of normal dual operator bimodules over 
von Neumann algebras (see [ER1]) extends to the non self-adjoint case.
\smallskip\noindent
We shall only use a little from Operator Space Theory. Anything undefined here can
be found in [Bl3], [BP], [ER2], or [Pi2]. 
We also refer the reader to [CS], [PS], and [CES] for some background on the Haagerup 
tensor product and operator bimodules, and to [BRS] for the abstract characterization
of unital operator algebras.
\bigskip\medskip\noindent
{\bf 2 - \underbar{A characterization of function modules.} }
\medskip\noindent
We first fix some tensor norm notation to be used in this section.
Given any two Banach spaces $X, Y$, we denote by
$Y \otimes X$ their algebraic tensor product. If $\alpha$ is a
tensor norm on $Y \otimes X$, the completion of $Y \otimes X$
under $\alpha$ is denoted by $Y\otimes_{\alpha} X$. The Banach space 
projective (resp. injective) tensor norm is denoted by $\alpha = \wedge$
(resp. $\alpha = \vee$). Then following [Pi1, Chapter 2], we  
denote by $\gamma_2$ the norm of factorization through Hilbert space.
Given any $1<q\leq\infty$, we let $g_q$ be the tensor norm defined 
in [DF, Section 12.7]. The key fact about this norm is not its definition
but its dual relationship with the notion of $p$-summing operators.
We recall that given a bounded operator $u\colon X\to Z$ between
Banach spaces, and a number $1\leq p <\infty$, $u$ is said to be $p$-summing
if there exists a constant $K>0$ such that for any finite family $(x_i)_i$ in $X$,
$\sum_i \Vert u(x_i)\Vert^p \leq K^p \sup\{\sum_i\vert\xi(x_i)\vert^p : \xi \in X^*, 
\Vert\xi\Vert\leq 1\}$. Moreover the infimum of all numbers $K$ for which this holds
is denoted by $\pi_p(u)$.
We refer e.g. to [Pi1, Chapter 1] for basic information on that notion.
Now the fundamental property of $g_q$ is the following.
Let $1<q\leq \infty$ and $1\leq p <\infty$ be such that $1/p + 1/q =1$.
Let $T\in \bigl(Y\otimes_{g_q} X\bigr)^*$ and let $u\colon X\to Y^*$ be the
corresponding linear mapping defined by $\langle u(x) , y\rangle = T(y\otimes x)$.
Then $u$ is $p$-summing and $\pi_p(u) = \Vert T\Vert$.
We finally recall the following contractive embedding, valid for any $1 < q\leq 2$:
$$
Y\otimes_{\wedge} X\subset Y\otimes_{g_q} X\subset
Y\otimes_{g_2} X\subset
Y\otimes_{\gamma_2} X\subset Y\otimes_{\vee} X\, .
\leqno(2.1)
$$
We shall also use the tensor norm $d_q$ (see [DF, Section 12.7]). 
The latter is defined by the isometric identity $Y\otimes_{g_q} X = X\otimes_{d_q} Y$.
\bigskip\noindent
Let $A$ be any unital Banach algebra and let $X$ be a left Banach $A$-module.
The module action gives rise to a linear mapping $m\colon A\otimes X\to X$ which
extends to a contraction from $A\otimes_{\wedge} X$ into $X$. If $\alpha$ is a tensor 
norm on $A\otimes X$ and if $m$ actually extends to a contraction
$A\otimes_{\alpha}X \to X$, we say that $X$ is a left $\alpha$-module over $A$.
When this holds with $\alpha=\vee$, we simply say that $X$ is a
left injective module over $A$.
Of course a similar definition can be given for right modules as well.
In [T], Tonge characterized function algebras as those unital Banach algebras $A$
for which the multiplication mapping on $A$ extends to a contraction
$A\otimes_{g_q}A\to A$, for some $1<q\leq \infty$. Using the main idea of the proof
of [T, Theorem 3], we shall give a similar representation result for
left $g_q$-modules. The relevant definition is the following.
\bigskip\noindent
{\bf Definition 2.1.} Let $X$ be a left Banach $A$-module. We say that
$X$ is a (left) {\it function module} over $A$, or a (left) function
$A$-module, if there is a compact space $K$, 
a linear isometry $i\colon X\to C(K)$, 
and a contractive unital homomorphism $\theta\colon A\to C(K)$, such that 
$i(a\cdotp x) =\theta(a)i(x)$ for any $a\in A,\ x\in X$. Right function modules
are defined similarly.
\bigskip\noindent
We warn the reader that our notion of function module does not coincide 
with the one defined and studied in [B2].
\smallskip\noindent
Before stating our result, we need to introduce a canonical mapping associated to any
Banach space.
Let $X$ be an arbitrary Banach space. Let ${\cal B}$ be the unit ball of $X^*$ endowed
with the $w^*$-topology (so that ${\cal B}$ is a compact space), and let 
$E_X\subset{\cal B}$ be the subset of its extreme points. Note that $\Vert\phi\Vert =1$
for any $\phi\in E_X$. We denote by $C_b(E_X)$ the 
unital commutative $C^*$-algebra of all bounded continuous functions on $E_X$. 
Then we let
$$
j_X\colon X\longrightarrow C_b(E_X)
$$
be the canonical mapping defined by $j_X(x)[\phi] =\langle \phi, x\rangle$ for any 
$x\in X, \phi\in E_X$. Then by the Hahn-Banach and Krein-Milman theorems, $j_X$
is an isometry.
\bigskip\noindent
{\bf Theorem 2.2.}
Let $A$ be a unital Banach algebra and let $X$ be a left Banach $A$-module.
The following four assertions are equivalent.
\smallskip
\item{(i)} $X$ is a left function module over $A$.
\smallskip
\item{(ii)} $X$ is a left injective module over $A$.
\smallskip
\item{(iii)} There exists $1< q\leq \infty$ such that $X$ is a
left $g_q$-module over $A$.
\smallskip
\item{(iv)} There exists a contractive unital
homomorphism $\theta\colon A\to C_b(E_X)$ such that:

$$
\forall x\in X,\, a\in A,\qquad  j_X(a\cdotp x) = \theta(a)j_X(x)\, .
\leqno(2.2)
$$
\bigskip\noindent
{\bf Proof.}
The unital Banach algebra $C_b(E_X)$ can be regarded as a $C(K)$-space, with
$K$ equal to the Stone-Cech compactification of $E_X$. Hence (iv) implies (i).
Clearly (ii) implies (iii), see e.g. (2.1). Now
assume (i), that is $A$ and $X$ satisfy Definition 2.1.
Let $a_1,\ldots, a_n$ and $x_1,\ldots,x_n$ be $n$-tuples 
in $A$ and $X$ respectively. The
multiplication mapping on $C(K)$ extends to a contraction on $C(K)\otimes_{\vee}
C(K)$. Furthermore, $\theta$ and $i$ are both contractions hence
$$
\Bigl\Vert \sum_{j=1}^n \theta(a_j)i(x_j)\Bigr\Vert_{C(K)}\leq
\Bigl\Vert \sum_{j=1}^n a_j\otimes x_j\Bigr\Vert_{A\otimes_{\vee}X}\, .
\leqno(2.3)
$$
Since $i$ is an isometry, it follows from the identity
$i(a\cdotp x) =\theta(a)i(x)$
that the left hand side of (2.3)
equals $\bigl\Vert m\bigl(\sum_j a_j\otimes x_j\bigr)\bigr\Vert$. This shows 
that $\bigl\Vert m\colon A\otimes_{\vee}X\to X\bigr\Vert\leq 1$, whence (ii).
\medskip\noindent
We now proceed to show that (iii) implies (iv), which is the main implication.
As mentioned before we use the main idea of Tonge's characterization of
function algebras ([T]), which in turn was inspired by an argument of Drury.
\smallskip\noindent
Let $1< q\leq \infty$ be such that
$\bigl\Vert m\colon A\otimes_{g_q}X\to X\bigr\Vert\leq 1$ and let 
$1\leq p<\infty$ be its conjugate number, that is $1/p+1/q =1$.
In the first part of this proof, we give ourselves a fixed
$\phi\in E_X$ and show that for any $a\in A$ and any 
$x,y\in X$, we have:
$$
\langle \phi, a\cdotp x\rangle\langle\phi,y\rangle =  
\langle \phi, a\cdotp y\rangle\langle\phi,x\rangle
\leqno(2.4)
$$
As before, we denote by ${\cal B}$ the 
unit ball of $X^*$. By our assumption, $m^*(\phi)$ is a contractive
functional on $A\otimes_{g_q} X$ hence the associated linear mapping
$u\colon X\to A^*$ defined by $\langle u(x), a\rangle = \langle\phi, a\cdotp x\rangle$
is $p$-summing, with $\pi_p(u)\leq 1$.
It therefore follows from the Pietsch factorization theorem (see e.g.
[Pi1, Theorem 1.3]) that there is a Radon probability measure $\mu$ on $\cal B$ such that
$$
\forall x\in X,\qquad \Vert u(x) \Vert^p \leq \int_{\cal B}\vert\langle x^*,x\rangle\vert^p\,
d\mu(x^*)\, .
\leqno(2.5)
$$
For any $a\in A$, we have
$\vert\langle\phi, a\cdotp x\rangle\vert\leq\Vert u(x)\Vert\Vert a\Vert$ 
hence we may deduce from (2.5) that there exists some $F_a\in L_q({\cal B};\mu)$
with norm less than $\Vert a\Vert$ such that
$$
\forall x\in X,\qquad
\langle\phi, a\cdotp x \rangle=\int_{\cal B} \langle x^*,x\rangle F_a(x^*)\, d\mu(x^*)\, .
\leqno(2.6)
$$
Let us simply denote by $F$ the function $F_a$ corresponding to the unit $a=1$.
Since $1\cdotp x = x$, we obtain
$$
\forall x\in X,\qquad
\langle\phi, x \rangle =\int_{\cal B} \langle x^*,x\rangle F(x^*)\, d\mu(x^*)\, ,
$$
with $\Vert F\Vert_{L_q({\cal B};\mu)}\leq 1$. Now using the fact that $\phi$
is an extreme point of ${\cal B}$, and following the argument in the proof of [T, Theorem 3],
we may deduce from the latter that 
$$
\forall x\in X,\qquad
\langle\phi, x \rangle = \langle x^*,x\rangle F(x^*)\quad \mu-\hbox{a.e.}\, .
\leqno(2.7)
$$
Let $a\in A$ and let $x,y\in X$. Then we can now check (2.4) as follows.
$$
\eqalign{
\langle \phi, a\cdotp x\rangle\langle\phi,y\rangle &
= \int_{\cal B} \langle\phi,y\rangle\langle x^*,x\rangle F_a(x^*)\, d\mu(x^*)
\quad\hbox{ by (2.6),}
\cr &
=\int_{\cal B} \langle x^*,y\rangle F(x^*)\langle x^*,x\rangle F_a(x^*)\, d\mu(x^*)
\quad\hbox{ by (2.7),}
\cr &
= \int_{\cal B} \langle\phi,x\rangle\langle x^*,y\rangle F_a(x^*)\, d\mu(x^*)
\quad\hbox{ by (2.7),}
\cr &
= \langle \phi, a\cdotp y\rangle\langle\phi,x\rangle
\quad\hbox{ by (2.6).}
\cr }
$$
For any $a\in A$, we may now define a function $\theta(a)\colon E_X\to\Cdb$ as follows.
Given any $\phi\in E_X$, we pick $x\in X$ so that $\langle\phi,x\rangle \not= 0$
(recall that $\phi\not= 0$) and we set:
$$
\theta(a)[\phi]=
{\langle \phi, a\cdotp x\rangle\over \langle\phi,x\rangle}\, .
$$
Indeed it follows from (2.4) that this definition does not depend on the choice of $x$.
Clearly $\theta(a)$ is a continuous function. Moreover since any $\phi\in E_X$ has norm $1$,
we have:
$$
\eqalign{
\bigl\vert  \theta(a)[\phi] \bigr\vert
& 
\leq\inf\{\vert \langle \phi, a\cdotp x\rangle\vert : 
x\in X, \,\vert \langle\phi,x\rangle\vert =1\}
\cr &
\leq \Vert a\Vert
\inf\{ \Vert x \Vert : x\in X, \,\vert \langle\phi,x\rangle\vert =1\} =\Vert a \Vert\, .
\cr }
$$
This shows that $\theta(a)\in C_b(E_X)$ for any $a\in A$ and that
$\theta\colon A\to C_b(E_X)$ is a contractive linear mapping. 
Moreover (2.4) can now be re-written as follows:
$$
\forall a\in A,\, x\in X, \phi\in E_X,\qquad
\theta(a)[\phi] \langle \phi, x\rangle = \langle \phi, a\cdotp x \rangle\, .
\leqno(2.8)
$$
If $a, b\in A$, we see by applying (2.8) three times 
that for any $\phi\in E_X$ and $x\in X$, we have:
$$
\eqalign{
\theta(ab)[\phi] \langle \phi, x\rangle  &
= \langle \phi, ab\cdotp x \rangle = \langle \phi, a\cdotp (b\cdotp x) \rangle
\cr &
= \theta(a)[\phi] \langle \phi, b\cdotp x\rangle
\cr &
= \theta(a)[\phi] \theta(b)[\phi] \langle \phi, x\rangle \, .
\cr }
$$
This shows that $\theta$ is a homomorphism. The latter is clearly unital.
Since (2.2) follows from (2.8), the proof is complete.
$\quad\diamondsuit$
\bigskip\noindent
{\bf Remark 2.3.} We may obviously write a version of Theorem 2.2 for
right modules. Namely, let $B$ be a unital Banach algebra and let $X$ be a
right Banach $B$-module. Then $X$ is a right function module
over $B$ if and only if $X$ is a right injective module over $B$
if and only if there exists $1< q\leq \infty$ such that $X$ is a
right $d_q$-module over $B$ if and only if 
there exists a contractive unital
homomorphism $\theta\colon B\to C_b(E_X)$ such that 
$\, j_X(x\cdotp b) = j_X(x)\theta(b)\,$ for any $x\in X, b\in B$.
\bigskip\noindent
As an application of the fact that we may use the same fixed
embedding $j_X$ in Theorem 2.2 and Remark 2.3, we get 
the following observation concerning bimodules.
\bigskip\noindent
{\bf Remark 2.4.} Let $A, B$ be two unital Banach algebras and assume 
that $X$ is both a left function $A$-module and a right function $B$-module.
Then there exist two contractive unital homomorphisms
$\theta_1\colon A\to C_b(E_X)$ and $\theta_2\colon B\to C_b(E_X)$ such that
$j_X(a\cdotp x)=\theta_1(a)j_X(x)$ and $j_X(x\cdotp b)= j_X(x) \theta_2(b)$
for any $a\in A,\, x\in X,\, b\in B$. We may then write
$$
j_X((a\cdotp x)\cdotp b) =\theta_1(a)j_X(x)\theta_2(b) = j_X(a\cdotp (x\cdotp b))\, .
$$
Since $j_X$ is 1-1, we deduce the associativity condition:
$$
\forall a\in A,\, x\in X,\, b\in B,\qquad
(a\cdotp x)\cdotp b = a\cdotp (x\cdotp b)\, .
$$
In other words the left $A$-module and right $B$-module actions on $X$ 
automatically induce an $A$-$B$-bimodule action on $X$.
\bigskip\noindent
{\bf Remark 2.5.}
\smallskip\noindent
{\bf (1)}$\,$
It is fairly easy to recover Tonge's characterization
of function algebras from our Theorem 2.2.
Indeed, let $1< q\leq \infty$ be a number and let $A$ be a unital Banach algebra.
Assume that the multiplication mapping $A\otimes A\to A$
extends to a contraction from $A\otimes_{g_q} A$ into $A$, and apply
Theorem 2.2 with $X=A$. Let $\theta\colon A\to C_b(E_A)$ be the contractive
unital homomorphism given by the latter. 
For any $a\in A$, we have $j_A(a)=\theta(a)j_A(1)$ hence
$$
\Vert a\Vert =\Vert j_A(a)\Vert = \Vert \theta(a)j_A(1)\Vert
\leq \Vert\theta(a)\Vert\Vert j_A(1)\Vert\leq\Vert\theta(a)\Vert\, .
$$
This shows that $\theta$ is actually an isometric unital homomorphism,
hence $A$ is a function algebra.
\medskip\noindent
{\bf (2)}$\,$ A ``non associative" version of Tonge's result
can be derived from the observation made in Remark 2.4.
Namely, let $A$ be a Banach space and assume that there is a bilinear
map $m \colon A\times A \to A$ and an element $e\in A$ such that $\Vert e\Vert =1$
and $m(a,e)=m(e,a)=a$ for any $a\in A$. Assume that for some $1<q_1,q_2\leq \infty$,
we both have 
$$
\bigl\Vert m\colon A\otimes_{g_{q_1}} A\to A\bigr\Vert\leq 1
\qquad\hbox{ and }\qquad
\bigl\Vert m\colon A\otimes_{d_{q_2}} A\to A\bigr\Vert\leq 1\, .
\leqno(2.9)
$$ 
Then $m$ is actually a multiplication on $A$, and $A$ equipped with $m$ is
a function algebra. 
\smallskip\noindent
Indeed, the proof of (2.8) above does not use the algebraic structure of $A$ hence 
under the assumption (2.9), it shows the existence of contractive linear maps
$\theta_1,\theta_2\colon A\to C_b(E_A)$ such that $\theta_1(a)j_A(b) =j_A(m(a,b))=
j_A(a)\theta_2(b)$ for any $a,b\in A$. Arguing as in Remark 2.4, we deduce the
associativity condition $m(m(a,c),b)=m(a,m(c,b))$. This should be compared
with [BRS, Corollary 2.4].
\bigskip\noindent
For any Banach space $X$, let us consider the multiplier algebra of $X$
defined by:
$$
{\rm Mult}(X) = \{ f\in C_b(E_X) : fj_X(X) \subset j_X(X)\}\, .
$$
Clearly ${\rm Mult}(X)\subset C_b(E_X)$ is a function algebra.
We refer e.g. to [B1, p. 54-55] for basic properties and an alternate definition
of this algebra. 
It follows from Theorem 2.2 that for any unital Banach algebra $A$,
there is a 1-1 correspondence between function $A$-modules actions on $X$ and contractive unital
homomorphisms from $A$ into ${\rm Mult}(X)$. Let us give a list of instructive examples
and observations, partly based on former results concerning multiplier algebras ([B1, B2, J]).
\bigskip\noindent
{\bf Examples 2.6.}
\smallskip\noindent
{\bf (1)}$\,$ We let $A$ be a unital Banach algebra.
If $\varphi\in A^*$ is a character of $A$, we may define a
$A$-module action on any $X$ by letting $a\cdotp x =\varphi(a)x$. Then $X$ is
an function $A$-module. Such a function module will be called {\it elementary}.
Conversely, if $\theta \colon A\to {\rm Mult}(X)$ is a contractive unital homomorphism and if
$\Psi$ is any character of ${\rm Mult}(X)$, then $\theta^*(\Psi)$
is a character of $A$. We deduce that the following are equivalent.
\smallskip
\item{(i)} $A$ admits characters.
\smallskip
\item{(ii)} There exists a non zero function $A$-module.
\smallskip
\item{(iii)} Any Banach space $X$ admits a function $A$-module action.
\medskip\noindent
{\bf (2)}$\,$ Let $X$ be a Banach space.
If ${\rm Mult}(X)=\Cdb$, then any function module action on
$X$ is elementary. This holds for instance if $X$ is strictly convex ([J, Corollary 3]),
or if $X$ is smooth ([B2, Proposition 5.2]), or if $X$ is an $L_1$-space 
([B2, p. 136]).
\medskip\noindent
{\bf (3)}$\,$ Let $X$ be a Banach space. If $X$ does not contain any isometric copy
of $c_0$, then ${\rm Mult}(X)$ is finite dimensional ([J]). This holds for instance if
$X$ is reflexive. Assuming this, let $n$ be the dimension of ${\rm Mult}(X)$.
Then ${\rm Mult}(X)=\ell_{\infty}^n$ hence
there exist subspaces $X_1,\ldots, X_n$ of $X$ such that
${\rm dim}\bigl({\rm Mult}(X_j)\bigr)=1$ for each $1\leq j\leq n$ and 
$X$ is equal to the $\ell_\infty$ direct sum $X=X_1\buildrel\infty\over\oplus
\cdots \buildrel\infty\over\oplus X_n$
([B2, Proposition 5.1]). This implies that if $A$ is a unital Banach algebra
admitting characters, the datum of a function $A$-module action on $X$
is equivalent to the datum of elementary actions of $A$ on each $X_j$.
More precisely, let $\varphi_1,\ldots,\varphi_n$ be characters on $A$.
We may define a function $A$-module action on $X$ by letting
$a\cdotp (x_1 +\cdots +x_n) = \varphi_1(a)x_1 +\cdots
+\varphi_n(a)x_n$ for any $a\in A,\ x_1\in X_1, \ldots ,x_n\in X_n$.
And conversely, any function $A$-module action on $X$ arises in that manner.
\medskip\noindent
{\bf (4)}$\,$
Let $X=B$ be a function algebra. Then ${\rm Mult}(B)=B$ ([B2, p. 123]). 
Thus for any unital Banach algebra $A$, there is a 1-1 correspondence between
$A$-module action on $B$ and contractive unital homomorphisms $\theta\colon A\to B$.
Given such a mapping $\theta$, the corresponding
action is given by $a\cdotp b = \theta(a)b$ for
$a\in A,\ b\in B$.
\medskip\noindent
{\bf (5)}$\,$
Let $A$ be the disc algebra and let $X$ be any Banach space.
Then the function $A$-module actions on $X$ are in an obvious correspondence
with elements of the closed unit ball of ${\rm Mult}(X)$.
\medskip\noindent
{\bf (6)}$\,$
The most interesting examples of function modules are perhaps those obtained
by the following process. Let $K$ be a compact space and let $A\subset C(K)$ be
a function algebra. Then for any set $S\subset C(K)$, the Banach space
$X=\overline{AS}=\overline{\rm Span}\{as : a\in A, s\in S\}$ is a left function
module over $A$.
\bigskip\noindent
{\bf Remark 2.7.}
If $A, B$ are unital Banach algebras, if $X$ is a left Banach $A$-module, 
and if $\rho\colon B\to A$ is a contractive unital homomorphism, then $X$ becomes a
left Banach $B$-module in a canonical way, namely $b\cdotp x =\rho(b)\cdot x$. 
We shall call this a {\it prolongation} of the $A$-module action.
Obviously, any prolongation of a function module is a function module.
In fact, for any $X$, the multiplication in $C_b(E_X)$ makes $X$ a left function
module over ${\rm Mult}(X)$, and every function module action 
on $X$ is a prolongation of this one.
\bigskip\noindent
We now turn to operator modules and operator bimodules in the sense of [CES].
Recall that the tensor norms $\gamma_2$, $d_2$, and $g_2$ are related to
the Haagerup tensor product of operator spaces as follows.
If $X$ and $Y$ are two Banach spaces, then the following isometric
equalities hold:
$$
MIN(Y)\otimes_h MIN(X) =Y\otimes_{\gamma_2} X\, ;
$$
$$
MAX(Y)\otimes_h MIN(X) =Y\otimes_{g_2} X\, ;\quad
MIN(Y)\otimes_h MAX(X) =Y\otimes_{d_2} X\, .
$$
See e.g. [Bl2, Theorem 3.1]. Applying Theorem 2.2 and Remarks 2.3 and 2.4, 
we therefore obtain the following.
\bigskip\noindent
{\bf Corollary 2.8.} Let $A$ and $B$ be two unital operator algebras, and let $X$ be a
Banach $A$-$B$-bimodule. Let $m\colon A\otimes X\otimes B\to X$, 
$m_1\colon A\otimes X\to X$ and $m_2\colon X\otimes B\to X$
be linearly defined by the module actions.
The following assertions are equivalent.
\smallskip
\item{(i)} There exist two numbers $1< q_1, q_2 \leq \infty$ such that
$$
\bigl\Vert m_1\colon A\otimes_{g_{q_1}} X\longrightarrow X\bigr\Vert\leq 1
\quad\hbox{ and }\quad
\bigl\Vert m_2\colon X\otimes_{d_{q_2}} B\longrightarrow X\bigr\Vert\leq 1\, .
$$
\smallskip
\item{(ii)} We have $\,\bigl\Vert m\colon MAX(A)\otimes_h MIN(X)\otimes_h MAX(B)\longrightarrow MIN(X)
\bigr\Vert_{cb}\leq 1$.
\smallskip
\item{(iii)} The operator space $MIN(X)$ is an operator $A$-$B$-bimodule.
\smallskip
\item{(iv)} We have $\,\bigl\Vert m_1\colon A\otimes_{\vee} X\longrightarrow X\bigr\Vert\leq 1\,$
and $\,\bigl\Vert m_2\colon X\otimes_{\vee} B\longrightarrow X\bigr\Vert\leq 1\,$.
\smallskip
\item{(v)} 
There exist two contractive unital homomorphisms
$\theta_1\colon A\to C_b(E_X),\, \theta_2\colon B\to C_b(E_X)$ such that

$$
\forall (a,x,b)\in A\times X\times B,
\qquad
j_X(a\cdotp x\cdotp b)=\theta_1(a)j_X(x)\theta_2(b)\, .
\leqno(2.10)
$$
\bigskip\noindent
Using part (1) in Example 2.6, we deduce an obstruction for a minimal operator space 
to be an operator module. 
\bigskip\noindent
{\bf Corollary 2.9.}
Let $A$ be an operator algebra without characters
(for instance, let $A=B(H)$, with dim$(H)\geq 2$).
Then no minimal operator space (except $(0)$) can be an operator $A$-module.
\bigskip\noindent
We now restrict to the case of function algebras, on which characters exist.
We obtain the following ``commutative" version of the Christensen-Effros-Sinclair
characterization of operator bimodules.
\bigskip\noindent
{\bf Corollary 2.10.} Let $A$ and $B$ be two function algebras, and let $X$ be a
Banach $A$-$B$-bimodule. 
The following assertions are equivalent.
\smallskip
\item{(i)} The operator space $MIN(X)$ is an operator $A$-$B$-bimodule.
\smallskip
\item{(ii)} There exist a unital commutative $C^*$-algebra $\cal C$, two
isometric unital homomorphisms $\pi_1\colon A\to {\cal C},\ 
\pi_2\colon B\to {\cal C}$, and a linear isometry $J\colon X\to {\cal C}$
such that $J(a\cdotp x\cdotp b) =\pi_1(a)J(x)\pi_2(b)$ for any
$a\in A, x\in X, b\in B$.
\bigskip\noindent
{\bf Proof.} Assume (i). By Corollary 2.8, 
there exist two contractive unital homomorphisms
$\theta_1\colon A\to C_b(E_X),\, \theta_2\colon B\to C_b(E_X)$ such that
(2.10) holds.
Since $A$ and $B$ are function algebras, we
may find unital commutative $C^*$-algebras $C_1, C_2$, and isometric unital
homomorphisms $\rho_1\colon A\to C_1,\, \rho_2\colon B\to C_2$. We also give ourselves
two characters $\varphi_1 \in A^*,\, \varphi_2\in B^*$.
\smallskip\noindent
We let ${\cal C} = C_1\buildrel\infty\over\oplus C_b(E_X)
\buildrel\infty\over\oplus C_2$ be the $\ell_{\infty}$-direct sum
of our three unital commutative $C^*$-algebras. Then ${\cal C}$ is a unital
commutative $C^*$-algebra as well. We now define $J\colon X\to {\cal C}$,
$\pi_1\colon A\to {\cal C}$ and $\pi_2\colon B\to {\cal C}$
by letting $J(x)=0\oplus j_X(x)\oplus 0$, 
$\pi_1(a)=\rho_1(a)\oplus\theta_1(a)\oplus\varphi_1(a)1$, and
$\pi_2(b)=\varphi_2(b)1\oplus\theta_2(b)\oplus\rho_2(b)$.
It is easy to check, using (2.10), that these mappings satisfy (ii).
$\quad\diamondsuit$
\bigskip\noindent
{\bf Remark 2.11.}
Clearly submodules of function modules, $L_\infty$-direct sums of function modules,
and prolongations of function modules (see Remark 2.7) are function modules.
\smallskip\noindent
Let $A$ be an operator algebra. The category ${ }_A FMOD$ of function $A$-modules, 
is a full subcategory of ${ }_A OMOD$, the operator $A$-modules. The relationship of the
class of function $A$-modules,
to the class ${ }_A OMOD$, seems to resemble the relationship of function algebras to
operator algebras. That is, unlike ${ }_A OMOD$, ${ }_A FMOD$ is not closed with respect to
certain elementary constructions. For example the quotient of a function module
by a closed submodule (which could perhaps be called a $Q$-module), is an operator $A$-module,
but not necessarily a function $A$-module. The module Haagerup
tensor product of two function modules (which one may view as an appropriate
quotient of the $\gamma_2$ tensor product), is not a function $A$-module (but is
an operator $A$-module). This also is connected with the fact that the category of 
minimal operator spaces is not closed under some basic operator space constructions.
\smallskip\noindent
Some of the ideas here also led to the work [Bl5]. In that paper singly generated
function modules are studied in a little bit more detail.
\bigskip\medskip\noindent
{\bf 3 - \underbar{A result on normal dual operator bimodules.} }
\medskip\noindent
The Christensen-Effros-Sinclair representation theorem for
operator $A$-$B$-bimodules (reviewed in Section 1 above)
is stated and proved in [CES] for unital $C^*$-algebras $A$ and $B$. 
However it is well-known that their result (namely [CES, Corollary 3.3]) 
extends with the same proof to the case when
$A$ and $B$ are merely non self-adjoint unital operator algebras.
In [ER1], a version of [CES, Corollary 3.3] is proved for the so-called
normal dual operator bimodules over a pair of von Neumann algebras $(A,B)$.
The aim of this section is to extend that result to the case when $A$ and $B$
are only assumed to be unital dual operator algebras. The new ingredient used to
treat this case is the $w^*$-version of Wittstock's representation theorem
for $w^*$-continuous completely bounded maps on dual operator algebras, recently
established in [L2].
\smallskip\noindent
Let us briefly recall or introduce a few
definitions. An operator space $X$ is called a dual operator
space if there is another operator space $Y$ such that $X=Y^*$ completely isometrically.
In that case there exists a $w^*$-continuous completely isometric embedding
of $X$ into some $B(H)$ (see [Bl3] or [ER1]). Accordingly, an operator algebra
$A$ is called a dual operator algebra if it can be represented algebraically
and completely isometrically as a $w^*$-closed subalgebra of $B(H)$ for some Hilbert
space $H$. We refer to [L1] for an abstract characterization of dual operator algebras.
\smallskip\noindent
Let $A$ and $B$ be two unital dual operator algebras and let $X$ be a dual operator space.
Assume that $X$ is an operator $A$-$B$-bimodule. Following [ER1], we say that
$X$ is a normal dual operator $A$-$B$-bimodule if in addition, the trilinear
mapping
$$
A\times X\times B\longrightarrow X, \quad (a,x,b)\mapsto a\cdotp x\cdotp b
$$
is separately $w^*$-continuous. The following is a representation theorem
for those bimodules. It extends
[ER1, Theorem 3.4] and [ER1, Theorem 4.1].
\bigskip\noindent
{\bf Theorem 3.1.}
Let $A$ and $B$ be two unital dual operator algebras and let $X$ be a normal
dual operator $A$-$B$-bimodule. 
\smallskip
\item{(1)} There exist a Hilbert space $H$, a $w^*$-continuous
complete isometry $J\colon X\to B(H)$, and $w^*$-continuous
completely isometric homomorphisms $\pi_1\colon A\to B(H),\
\pi_2\colon B\to B(H)$ such that $J(a\cdotp x\cdotp b)
=\pi_1(a)J(x)\pi_2(b)$ for any $a\in A$, $x\in X$, $b\in B$.
\smallskip
\item{(2)} If $A=B$, we may find a $w^*$-continuous and
completely isometric unital homomorphism $\pi\colon A\to B(H)$
such that (1) holds with $\pi_1=\pi_2=\pi$.
\bigskip\noindent
{\bf Proof.} Let $m\colon A\otimes_h X\otimes_h B\to X$ be the completely contractive
mapping induced by the bimodule action. We let $j\colon X\to B(K)$ be a 
$w^*$-continuous completely isometric embedding, for some Hilbert space $K$. Then 
$jm$ is completely contractive hence by the factorization theorem of multilinear completely
bounded maps ([CS, PS]), there exist two Hilbert spaces $E, F$, and
three completely contractive maps $\alpha\colon A\to B(E,K)$,
$v \colon X\to B(F,E)$, and $\beta\colon B\to B(K,F)$, such that:
$$
\forall a\in A,\, x\in X,\, b\in B,\qquad
j(a\cdotp x\cdotp b) = \alpha(a) v(x)\beta(b)\, .
\leqno(3.1)
$$
A key fact is that since $m$ is assumed to be separately $w^*$-continuous
and $j$ is $w^*$-continuous, it follows from [ER3, Theorem 3.1] that
our three maps $\alpha, v, \beta$ can be chosen to be $w^*$-continuous.
By [L2, Proposition 3.4], we can then factor $\alpha$ and $\beta$ as follows.
There exist a Hilbert space $H_1$, a $w^*$-continuous unital homomorphism
$\theta_1\colon A\to B(H_1)$ and two contractions $E\buildrel S\over\longrightarrow H_1
\buildrel R\over\longrightarrow K$ such that $\alpha(a)=R\theta_1(a)S$. Likewise,
there exist a Hilbert space $H_2$, a $w^*$-continuous unital homomorphism
$\theta_2\colon B\to B(H_2)$ and two contractions $K\buildrel U\over\longrightarrow H_2
\buildrel T\over\longrightarrow F$ such that $\beta(b)=T\theta_2(b)U$.
Moreover, changing $H_1$ and $H_2$ to subspaces if necessary, we may assume
(as in the proof of [L2, Proposition 3.4]) that:
$$
H_1 =\overline{\rm Span}\bigl\{\theta_1(A)^*R^*(K)\bigr\}
\quad\hbox{ and }\quad
H_2 =\overline{\rm Span}\bigl\{\theta_2(B)U(K)\bigr\}\, .
\leqno(3.2)
$$
Let us now define $u\colon X\to B(H_2,H_1)$ by letting $u(x)=Sv(x)T$. By
construction, $u$ is $w^*$-continuous and completely contractive. Moreover
we see from (3.1) that:
$$
\forall a\in A,\, x\in X,\, b\in B,\qquad
j(a\cdotp x\cdotp b) = R\theta_1(a) u(x)\theta_2(b)U\, .
\leqno(3.3)
$$
Applying (3.3) twice it is fairly easy to see that for any
$a', a\in A,\, x\in X,\, b,b'\in B$, we have:
$$
R\theta_1(a')\theta_1(a) u(x)\theta_2(b)\theta_2(b')U =
R\theta_1(a') u(a\cdotp x\cdotp b)\theta_2(b')U\, .
$$
It therefore follows from (3.2) that for any 
$a\in A,\, x\in X,\, b\in B$, we have:
$$
u(a\cdotp x\cdotp b) = \theta_1(a) u(x)\theta_2(b)\, .
\leqno(3.4)
$$
Let $H=H_1\buildrel{2}\over\oplus H_2$ be the Hilbert direct sum of $H_1$
and $H_2$, and let $\pi_1\colon A\to B(H)$,
$\pi_2\colon B\to B(H)$, and $J\colon X\to B(H)$ be defined by letting:
$$
J(x)=\pmatrix{ 0 & u(x)\cr 0 & 0 \cr }\, ,\quad
\pi_1(a)=\pmatrix{ \theta_1(a) & 0\cr 0 & 0 \cr }\, ,\quad
\pi_2(b)=\pmatrix{ 0 & 0\cr 0 & \theta_2(b) \cr }\, .
$$
Clearly $\pi_1$ and $\pi_2$ are $w^*$-continuous completely
contractive homomorphisms. Moreover the identity (3.4) implies that:
$$
\forall a\in A,\, x\in X,\, b\in B,\qquad
J(a\cdotp x\cdotp b) = \pi_1(a) J(x)\pi_2(b)\, .
\leqno(3.5)
$$
Applying (3.3) with $a=1$ and $b=1$, we see that $u$ is actually a complete
isometry hence $J$ is a $w^*$-continuous complete isometry. 
The assertion (1) is thus proved, up to the fact that $\pi_1$ and $\pi_2$
are merely completely contractive instead of being completely isometric.
To get this completely isometric condition, it suffices to apply a direct
sum argument as e.g. in the proof of Corollary 2.10 above.
\medskip\noindent
We now turn to the proof of (2). Assume that $A=B$ and that $\theta_1$,
$\theta_2$ and $u$ are constructed as above. Again we let
$H=H_1\buildrel{2}\over\oplus H_2$ and this time, we define
$\pi\colon A\to B(H)$ and $J\colon X\to B(H)$ by letting:
$$
J(x)=\pmatrix{ 0 & u(x) \cr 0 & 0 \cr }\, ,\quad
\pi(a)=\pmatrix{ \theta_1(a) & 0\cr 0 & \theta_2(a) \cr }\, .
$$
Then $\pi$ is a $w^*$-continuous completely contractive unital
homomorphism and (3.5) holds with $\pi_1=\pi_2=\pi$. We can thus
conclude as in the proof of (1).
$\quad\diamondsuit$
\bigskip\noindent
{\bf Remark 3.2.} In [Bl4, Theorem 2.2], the first author gave a
completely isomorphic version of 
the Christensen-Effros-Sinclair representation theorem, as follows.
Let ($A, B$) be a pair of (possibly non unital) operator
algebras, and let $X$ be an operator space. Assume that $A$ is an $A$-$B$-bimodule
and that the bimodule action is completely bounded from $A\otimes_h X \otimes_h B$
into $X$.
Then one can find complete isomorphisms $J\colon X\to B(H)$,
$\pi_1\colon A\to B(H)$ and $\pi_2\colon B\to B(H)$ such that $\pi_1,
\pi_2$ are homomorphisms and 
$J(a\cdotp x\cdotp b)
=\pi_1(a)J(x)\pi_2(b)$ for any $a\in A$, $x\in X$, $b\in B$.
Arguing as in the proof of [Bl4, Theorem 2.2] and using the completely
isomorphic characterization of dual operator algebras given in [L1, Theorem 7],
one can show that if $X$ is a dual operator space, if $A$ and $B$ are dual operator
algebras, and if the bimodule action is separately $w^*$-continuous, then that
representation can be achieved with the extra condition that $\pi_1$, $\pi_2$, and
$J$ are $w^*$-continuous. We leave the details to the reader.
\bigskip\noindent
{\bf Remark 3.3.} Comparing our results in Section 2 and 3, it is tempting to try
to study an appropriate notion of dual function module. The first question arising here
is the following apparently open problem. Let $A$ be a function algebra and assume that
$A$ is a dual space. Does there exist 
a $w^*$-continuous isometric unital homomorphism of $A$ into some $L_\infty$-space?
\vskip 1cm
\noindent
{\bf Acknowledgments.} This paper was written while the second author was visiting the 
University of Houston. He wishes to thank the Department of Mathematics of UH for its 
warm hospitality.
\vskip 1cm
\centerline{\bf REFERENCES}
\bigskip
\baselineskip=11pt
\item{[B1]} E. Behrends, M-structure and the Banach-Stone theorem, Springer Lecture Notes
736, Berlin-Heidelberg-New-York, 1979.
\smallskip
\item{[B2]} E. Behrends, Multiplier representations and an application to the problem
whether $A\otimes_{\varepsilon} X$ determines $A$ and/or $X$, Math. Scand. 52 (1983), 117-144.
\smallskip
\item{[Bl1]} D.P. Blecher, Commutativity in operator spaces, Proc. 
Amer. Math. Soc 109 (1990), 709-715.
\smallskip
\item{[Bl2]} D.P. Blecher, Tensor products of operator spaces II, Canad. J. Math. 44
(1992), 75-90.
\smallskip
\item{[Bl3]} D.P. Blecher, The standard dual of an operator space, Pacific J. Math.
153 (1992), 15-30.
\smallskip
\item{[Bl4]} D.P. Blecher, A generalization of Hilbert modules, J. Funct. Anal.
136 (1996), 365-421.
\smallskip
\item{[Bl5]} D.P. Blecher, The Shilov boundary of an operator space - and
applications to the characterization theorems and Hilbert $C^*$-modules, Preprint (1999).
\smallskip
\item{[BP]} D.P. Blecher and V.I. Paulsen, Tensor products of operator spaces,
J. Funct. Anal. 99 (1991), 262-292.
\smallskip
\item{[BRS]} D.P. Blecher, Z.-J. Ruan and A.M. Sinclair, A characterization of operator algebras,
J. Funct. Anal. 89 (1990), 188-201.
\smallskip
\item{[CES]} E. Christensen, E.G. Effros and A.M. Sinclair, Completely bounded multilinear maps and
$C^*$-algebraic cohomology, Invent. Math. 90 (1987), 279-296.
\smallskip
\item{[CS]} E. Christensen and A.M. Sinclair, Representation of completely bounded multilinear operators,
J. Funct. Anal. 72 (1987), 151-181.
\smallskip
\item{[DF]} A. Defant and K. Floret, Tensor norms and operator ideals, North-Holland, Amsterdam, 1993.
\smallskip
\item{[ER1]} E.G. Effros and Z.-J. Ruan, Representations of operator modules and their
applications, J. Operator Theory  19 (1988), 137-157.
\smallskip
\item{[ER2]} E.G. Effros and Z.-J. Ruan, A new approach to operator spaces,
Canadian Math. Bull. 34 (1991), 329-337.
\smallskip
\item{[ER3]} E.G. Effros and Z.-J. Ruan, Operator convolution algebras: an approach to
quantum groups, Unpublished (1991).
\smallskip
\item{[J]} K. Jarosz, Multipliers in complex Banach spaces and structure of the
unit balls, Studia Math. T. LXXXVII (1987), 197-213.
\smallskip
\item{[L1]} C. Le Merdy, An operator space characterization of dual operator algebras,
Amer. J. Math. 121 (1999), 55-63.
\smallskip
\item{[L2]} C. Le Merdy, Finite rank approximation and semidiscreteness for linear operators,
Annales Inst. Fourier, to appear.
\smallskip
\item{[PS]} V.I. Paulsen and R.R. Smith, Multilinear maps and tensor norms on operator systems,
J. Funct. Anal. 73 (1987), 258-276.
\smallskip
\item{[Pi1]} G. Pisier, Factorization of linear operators and geometry of Banach spaces,
CBMS Series 60 (Amer. Math. Soc., Providence, R.I.), 1986.
\smallskip
\item{[Pi2]} G. Pisier, An introduction to the theory of operator spaces,
Preprint (1997).
\smallskip
\item{[T]} A.M. Tonge, Banach algebras and absolutely summing operators, Math. Proc.
Cambridge Philos. Soc. 80 (1976), 465-473.
\vskip 1cm
\noindent
{\bf David BLECHER:} {\it Department of Mathematics, University of Houston, Houston TX 77204-3476.}
\smallskip\noindent
{\bf Christian LE MERDY:} {\it D\'epartement de math\'ematiques, Universit\'e de Franche-Comt\'e,
25030 Besan\c con Cedex (France).}

\bye